\def\sqr#1#2{{\vcenter{\hrule height.#2pt              
     \hbox{\vrule width.#2pt height#1pt\kern#1pt
     \vrule width.#2pt}
     \hrule height.#2pt}}}
\def\square{\mathchoice\sqr{5.5}4\sqr{5.0}4\sqr{4.8}3\sqr{4.8}3}
\def\qed{\hskip4pt plus1fill\ $\square$\par\medbreak}

\input epsf.sty

\magnification\magstep1\overfullrule=0pt

\centerline{\bf Linear Recurrences in the Degree}

\centerline{\bf Sequences of Monomial  Mappings}

\bigskip
\centerline{ Eric Bedford and Kyounghee Kim}
\bigskip
\noindent{\bf \S1.  Introduction.  }  Let $A$ denote a $k\times k$ matrix of rank $k$ which has integer entries.  The monomial map  $f_A:{\bf C}^k\to {\bf C}^k$ defined by 
$$f_A(x)=x^A=\left( \prod_j x_j^{a_{1,j}}, \dots, \prod_j x_j^{a_{n,j}}\right)\eqno(1.1)$$
is a dominant rational map.  The iterates are given by $f_A^n=f_A\circ\cdots\circ f_A=f_{A^n}$.   If  $A\in GL(n,{\bf Z})$, then $f_A$ is birational, and   $f_A^{-1}=f_{A^{-1}}$.  A rational map $f$ on projective space ${\bf P}^k$ induces a linear map $f^*$ on $H^{p,p}({\bf P}^k;{\bf Z})\cong {\bf Z}$.  We define the degree of $f^n$ in codimension $p$ to be $d^{(n)}_p:=(f^n)^*|H^{p,p}({\bf P}^k)$; or equivalently (see [RS]), if $\omega$ is a K\"ahler form on ${\bf P}^k$ normalized so that $\int_{{\bf P}^k}\omega^k=1$, then
$$d^{(n)}_p=\int_{{\bf P}^k} f^*(\omega^p)\wedge\omega^{n-p}.$$

If $f:{\bf P}^k\to{\bf P}^k$ is rational, if $X$ is a compact K\"ahler manifold, and if $\pi:X\to{\bf P}^k$ is holomorphic and bimeromorphic, then we have a map $\tilde f=\pi^{-1}\circ f\circ \pi:X\to X$.  There will be an induced linear map
$$\tilde f^*:H^{p,p}(X)\to H^{p,p}(X).$$
Let $\chi (x)=x^m+\alpha_{m-1}x^{m-1}+\cdots+\alpha_0\in{\bf Z}[x]$ be the characteristic polynomial of $\tilde f^*|H^{p,p}$.  If we have
$$(\tilde f^n)^*=(\tilde f^*)^n\ \ {\rm on\ \ }H^{p,p}(X)\eqno(1.2)$$
then by [DF, Corollary 2.2]  the sequence $\{d^{(n)}_p\}_{n\in{\bf N}}$ satisfies the linear recurrence
$$ d_{p}^{(n+m)} + \alpha_{m-1}d^{(n+m-1)}_{p}+\cdots+\alpha_0 d^{(n)}_p=0  \eqno(1.3)$$
for all $n\in{\bf N}$.

In dimension $k=2$, there is only the case $p=1$ to consider.  Favre [F] has given necessary and sufficient conditions for a monomial map  in dimension 2  to have a regularization $\pi:X\to{\bf P}^2$ satisfying (1.2).  Diller and Favre [DF] showed that for every bimeromorphic surface map there is such  a regularization.  Favre and Jonsson [FJ] have shown that the degree sequence of a polynomial map of ${\bf C}^2$  always satisfies (1.3).

In dimension $k=3$,   Hasselblatt and Propp [HP] showed that there is  a matrix  $A\in GL(3,{\bf Z})$ such that the degree sequence $\{d^{(n)}_1\}=\{{\rm deg}(f_A^n)\}$ does not satisfy any linear recurrence. 

In homogeneous coordinates we have
$$f_A[x_0:x_1:\dots:x_k] = [1:\prod_j (x_j/x_0)^{a_{1,j}}: \cdots: \prod_j (x_j/x_0)^{a_{n,j}} ].\eqno(1.4)$$
and if we rewrite $f_A$ so that the coordinates are homogeneous polynomials, then their degree is $d^{(1)}_1=D(A)$, with 
$$D(A):= \max(0,\sum_{j=1}^k a_{1,j},\dots,\sum_{j=1}^k a_{n,j}) + \sum_{j=1}^k \max(0,-a_{1,j},\dots,-a_{n,j}).\eqno(1.5)$$
It is evident that there is a set ${\cal C}$ consisting of $(k+1)^{k+1}$ linear functionals $L_C:{\bf M}_k\to{\bf R}$ on the space of  real $k\times k$ matrices such that 
$$D(A)=\max\{L_C(A):C\in{\cal C}\}\eqno(1.6)$$
For any linear functional  $L$ on the set of $k\times k$ real matrices,  the sequence $\{L(A^n)\}_{n\in{\bf N}}$ satisfies the linear recurrence (1.3)  (see Lemma 2.1), where $\chi(x)$ is the characteristic polynomial of $A$.  It follows that the degree sequence $\{d^{(n)}_1\}=\{D(A^n)\}$ is ``almost'' the solution to a linear recurrence: it is the finite maximum over $C\in {\cal C}$ of the sequences $\{L_C(A^n)\}_{n\in{\bf N}}$, each of which satisfies (1.3) but which may have different initial conditions.  Another way of describing this phenomenon is that the space of matrices is divided into different cells defined by 
$$S_C:=\{M\in {\bf M}_k: D(M)=L_C(M)\}.\eqno(1.7)$$  For each $n$, $A^n$ belongs to one of the cells $S_{C(n)}$, so $d^{(n)}_1=L_{C(n)}(A^n)$.  Although there are only finitely many cells,  the cell $C(n)$ might change in a sufficiently irregular way that $\{d^{(n)}_1\}_{n\in{\bf N}}$ does not satisfy any linear recurrence at all.  This is the approach taken in [HP] and also used here in \S3.

Here we extend the results of [HP]  to obtain the following:
\proclaim Theorem 1.1.  Let $A$ be a $k\times k$ integer matrix with rank $k$.  Suppose that  for every eigenvalue $\lambda$ of $A$ with  $|\lambda|\ge1$, $\bar\lambda/\lambda$ is a root of unity.  Then $\{d^{(n)}_1\}_{n\in{\bf N}}$ satisfies a linear recurrence.  On the other hand, suppose that  the eigenvalues of largest modulus consist of a conjugate pair $\lambda$ and $\bar\lambda$ of simple eigenvalues, and that $\bar\lambda/\lambda$ is not a root of unity.  Then $\{d^{(n)}_1\}_{n\in{\bf N}}$ does not satisfy any linear recurrence.

By the duality between $H^{1,1}$ and $H^{k-1,k-1}$ we see that $f_A^*|H^{k-1,k-1}$ is dual to $(f_A)_*|H^{1,1} = f_{A^{-1}}^*|H^{1,1}$.  Thus in dimension 3, $\{d_2(f_A^n)\}= \{d_1(f_{A^{-1}}^n)\}$.   By this duality,  Theorem 1.1 gives a rather complete treatment of the cases in dimension 3:
\proclaim Theorem 1.2.  Suppose that the matrix $A\in GL(3,{\bf Z})$ has no eigenvalues of  modulus one.  In case all three eigenvalues are real, then both $\{d^{(n)}_1\}_{n\in{\bf N}}$ and $\{d^{(n)}_2\}_{n\in{\bf N}}$ satisfy linear recurrences.  In case there is a non-real eigenvalue,  we may write the eigenvalues as $\lambda$, $\bar\lambda$, and $\pm|\lambda|^{-2}$.  If $|\lambda|>1$, then $\{d^{(n)}_1\}_{n\in{\bf N}}$ does not satisfy a linear recurrence, but $\{d^{(n)}_2\}_{n\in{\bf N}}$ does.  And vice versa if  $|\lambda|<1$.

In dimension 4 Theorem~1.1 and duality, applied to the map $f_A:{\bf P}^4\to{\bf P}^4$, give:
\proclaim Theorem 1.3.  Suppose that the eigenvalues of  $A\in GL(4,{\bf Z})$ are two conjugate pairs $\lambda_j,\bar\lambda_j$,  $j=1,2$, such that  neither $\bar\lambda_1/\lambda_1$ nor $\bar\lambda_2/\lambda_2$  is a root of unity,  and $|\lambda_1|<1<|\lambda_2|$.  Then neither $\{d^{(n)}_1\}_{n\in{\bf N}}$ nor  $\{d^{(n)}_3\}_{n\in{\bf N}}$ is given by a linear recurrence.  In particular, there is no map $\tilde f_A$ satisfying (1.2) for $p=1$ or $p=3$.

In \S2 we prove the first part of Theorem~1.1, which involves matrices with essentially real eigenvalues; the second part is proved in \S3.

\bigskip\noindent{\bf \S2. Monomial mappings:  linear recurrence for $d^{(n)}_1$.}   Recall the linear functions $L_C$ in (1.6) and the cells $S_C$ in (1.7).  From the form of $D(A)$ in (1.5) we see that the boundary $\partial S_C$ is contained in a union of hyperplanes in the set $\{A=(a_{i,j})\}\subset {\bf M}_k$ of $k\times k$ matrices.  These hyperplanes are defined, for fixed $i$, $j$, and $m$, by   $\{a_{i,j}=0\}, \{a_{i,j} = a_{m,j}\}, \{\sum_\sigma a_{i,\sigma}=0\},$ and $\{\sum_\sigma a_{i,\sigma} = \sum_\sigma a_{m,\sigma}\}$.   We now make two observations about linear recurrences.
\proclaim Lemma 2.1.  Let a $k\times k$ matrix $A$ be given.  For fixed $1\le i,j\le k$, the sequence of  $i,j$ elements $\{(A^n)_{i,j}\}_{n\in{\bf N}}$ satisfies (1.3), and thus for fixed  $C\in{\cal C}$ the sequence $\{L_C(A^n)\}_{n\in{\bf N}}$ satisfies (1.3)

\noindent{\it Proof. }  By the Cayley Hamilton Theorem, the characteristic polynomial satisfies $\chi(A)=0$.  Thus each entry $(A^n)_{i,j}$ satisfies (1.3), and so the Lemma follows.  \qed

\proclaim Corollary 2.2. Suppose $A$ is a $k\times k$ integer matrix, and $N<\infty$ is such that  $A^n \in S_C$ for some fixed $C\in{\cal C}$ and for all $n\ge  N$. Then the degree sequence $\{d_1^{(n)}\}$ for $f_A$ satisfies a linear recurrence relation with constant coefficients. 

Let $\lambda_1, \dots, \lambda_t$ denote the distinct non-zero eigenvalues of A, and let $\mu_j\ge1$ denote the size of  the largest Jordan block corresponding to the eigenvalue $\lambda_j$.  Then there exist  constants $\alpha_{i,j}(s,\ell)$ which do not depend on $n$ such that:
$$(A^n)_{i,j} = \sum_{\ell=1}^t \left( \sum_{s=0}^{\mu_\ell-1} \alpha_{i,j}(s,\ell) {n \choose s} \lambda_\ell^{n-s} \right)\eqno{(2.1)}$$

\noindent{\bf Case 1: All eigenvalues are positive: $\lambda_1 >  \cdots> \lambda_t >0$.}  For a linear functional $L:{\bf M}_k\to{\bf R}$, we define $Q_L(n):=L(A^n)$.  With reference to (2.1) we set 
$$Q_{L,\ell}(n)= \sum_{s=0}^{\mu_\ell-1} L(\alpha_{i,j}(s,\ell) ){n \choose s} \lambda_\ell^{n-s},$$
and so we have
$$Q_L(n)=   \sum_{\ell=1}^t Q_{L,\ell}(n).\eqno{(2.2)}$$      
\proclaim Lemma 2.3. In Case 1, there is an integer $N$ such that either (i) $Q_L(n) \ge 0 $ for all $n \ge N$ or (ii) $Q_L(n) \le 0$ for all $n \ge N$.

\noindent{\it Proof.}  Let us look at the form of $Q_{L,\ell}(n)$.  It is the sum of terms $ {n\choose s}$ which are polynomial in $n$ and $\lambda_\ell^{n-s}$, which are exponential in $n$.  Thus $Q_{L,\ell}(n)$ may be written as a polynomial $P_\ell(n)$ in $n$ multiplied by $\lambda_\ell^{n}$.  Let $\ell_0$ denote the first value of $\ell$ for which $P_{\ell_0}$ is not the zero polynomial.  It follows that $|P_{\ell_0}(n)|$ grows like a power of $n$, and  times $\lambda_{\ell_0}^n$.  Summing ver the remaining $\ell$, we see that $Q_L(n)$ has the same growth.  Thus $Q_L(n)$ is ultimately $\ge0$ or $\le0$.  \qed

\proclaim Lemma 2.4. In Case 1, there exists a positive integer $N$ such that for all $n \ge N$, $A^n$ belongs to one particular cell.

\noindent{\it Proof.}  For $C\in{\cal C}$, we consider a linear functional that defines one of the sides of a cell $S_C$.  By Lemma 2.3, we know that the sequence $L(A^n)$ is ultimately $\ge0$ or $\le 0$.  We apply this to all of the linear functionals defining $S_C$, and we find that for all $n$ sufficiently large,  either $A^n\in S_C$ or $A^n\notin S_C$.  On the other hand, the set of all $S_C$, $C\in{\cal C}$, exhausts ${\bf M}_k$, so $A^n$ must belong  to one particular cell.\qed 

\noindent{\bf Case  2:  Eigenvalues whose modulus is $\ge1$  are positive.}  Let us choose $t_0$ such that the eigenvalues of $A$ are given by  $\lambda_1 > \lambda_2 > \cdots > \lambda_{t_0} \ge 1>| \lambda_{t_0+1}|\ge  \cdots \ge | \lambda_{t_1}| >0 $.   As before, we let $L:{\bf M}_k\to{\bf R}$ be a linear functional, but now we suppose that it has integer coefficients.  Then with $Q_L$ and $Q_{L,\ell}$ as above, we define $R_L(n)$ by
$$Q_L(n) = \sum_{\ell =1}^{t_0}Q_{L,\ell}(n) + R_L(n).$$

\proclaim Lemma 2.5.  In Case 2, there is an integer $N$ such that either (i) $Q_L(n) \ge 0 $ for all $n \ge N$ or (ii) $Q_L(n) \le 0$ for all $n \ge N$.

\noindent{\it Proof.}  As in the Proof of Lemma 2.3, each $Q_{L,\ell}(n)$ is either identically zero or grows like a power of $n$ times $\lambda_\ell^n$.   The  conclusion of the Lemma must hold, then, unless $Q_{L,\ell}(n)=0$ for $\ell$ equal to $1,\dots,t_0$.  This means that $Q_L(n)=R_L(n)$.  But now we recall that $L$ has integer coefficients, so that $Q_L(n)\in{\bf Z}$.  On the other hand, $R_L(n)$ is made with exponentials with modulus $<1$, so we have $R_L(n)\to 0$ as $n\to\infty$.  We this can happen only if $Q_L(n)=0$ for all $n$.  \qed

We observe that the hyperplanes bounding $S_C$ are defined by functions with integer coefficients, so as in Lemma 2.4 we have:
\proclaim Lemma 2.6.  In Case 2, there exists a positive integer $N$ such that for all $n \ge N$, $A^n$ belongs to one particular cell.

\proclaim Theorem 2.7.     Let $A$ be a $k\times k$ integer matrix with rank $k$.  Suppose that  for every eigenvalue $\lambda$ of $A$ with  $|\lambda|\ge1$, $\bar\lambda/\lambda$ is a root of unity.  Then $\{d^{(n)}_1\}_{n\in{\bf N}}$ satisfies a linear recurrence.   If the eigenvalues $|\lambda|\ge1$ are all positive, then the degree sequence satisfies the linear recurrence relation given by the characteristic polynomial of $A$.

\noindent{\it Proof.}    For each such eigenvalue $\lambda$ with $|\lambda|\ge1$ we choose $\tau>1$ such that $(\lambda/|\lambda|)^\tau=1$.  Let $\tilde \tau$ be the least common multiple of all such $\tau$.  Now the eigenvalues of $A^\tau$ are in Case 2.  By Lemma 2.6, we find that for fixed $\nu$ with $0\le \nu\le\tilde\tau-1$, $\{ d_1^{\tilde\tau n+\nu} )\}$ satisfies a linear recurrence whose coefficients are given by the characteristic polynomial of $A^\tau$.  Thus the full sequence also satisfies a linear recurrence.  \qed

\bigskip\noindent{\bf \S3. Monomial mappings:  no linear recurrence for $d^{(n)}_1$.}   We will use a fact from Combinatorics ( see [S, Chapter 4]):  If $\{c_k\}$ and $\{d_k\}$ both satisfy  linear recurrence relations, then the indices $k$ for which $c_k-d_k = 0$ is eventually periodic. 

\proclaim Proposition 3.1. Suppose that $A= (a_{i,j})$ is an integer matrix of rank $k$, and suppose $A$ has exactly two eigenvalues $\delta,\bar\delta$ of maximum modulus, and $\bar\delta/\delta$ is not a root of unity. Then the degree sequence $d_1^{(n)}$ for $f_A$ does not satisfy a linear recurrence. 

\noindent{\it Proof.}  Let $m$ denote the size of the largest Jordan block with eigenvalue $\delta$ or $\bar\delta$.  We choose $C\in{\cal C}$ such that $A^n\in{\cal C}$ for infinitely many $n$.  Writing $\lambda_\ell$, $1\le \ell\le t$ for the other eigenvalues of $A$, we have
$$c_n:=L_C(A^n) =  \sum_{s=0} ^{m-1} {n \choose s} Re \left( \beta(s) \delta^{n-s}\right) + \sum_{\ell=1}^t \left( \sum_{s=0}^{\mu_\ell-1} \alpha_{i,j}(s,\ell) {n \choose s} \lambda_\ell^{n-s} \right) $$
For the values of $n$ such that $A^n\in S_C$ we have $d^{(n)}_1=L_C(A^n)$.  Since $d^{(n)}_1$ grows like $|\delta|^n$, and since $|\lambda_j|<|\delta|$, we see that not all the coefficients $\beta(s)$ can be equal to zero.

By Lemma 2.1, $\{c_n\}$ satisfies a linear recurrence.  If $\{d_1^{(n)}\}$ also satisfies a (possibly different) linear recurrence, then by the combinatorial fact above, the indices for which $d^{(n)}_1=c_n$ are eventually periodic.  This means that they agree for $n$ belonging to an arithmetic progression $An+B$.  Now we write $\eta=(\delta/|\delta|)^A$.  Since $\eta$ is not a root of unity, the numbers $Re(\eta^n)$ are dense in the unit circle.  Restricting to this arithmetic sequence, we have
$$d^{(An+B)}_1=c_{An+B} =   |\delta|^{An}\sum_{s=0} ^{m-1} {An + B \choose s} Re \left(\delta^{B-sA} \beta(s) \eta^{n}\right) + O(|\lambda|^{An}).$$
We let $s_0$ denote the largest value of $s$ such that $\beta(s_0)\ne0$, and for large $n$ this will give the dominant term in the summation.  However, there are arbitrarily large values of $n$ for  which  
$Re(\cdots)\le -{1\over 2} |\delta^{B-s_0A}\beta(s_0)|$.  So, for such a value of $n$ which is sufficiently large, $d_1^{(n)}$ will be negative, which is a contradiction.  \qed

\noindent{\it Proof of Theorem 1.1.}  This follows directly from Theorems 2.7 and 3.1.
\medskip
\noindent{\it Proof of Theorem 1.2.}  The characteristic polynomial of $A$ is an irreducible cubic so $\bar\lambda/\lambda$ cannot be a root of unity, and Theorem 1.2 follows from Theorem 1.1.
\medskip
\noindent {\bf Example of Hasselblatt and Propp.}  The example given in [HP]  is
$$f_A (x_1,x_2,x_3) := ( {x_2 \over x_1},{x_3 \over x_1},x_1).$$ 
The eigenvalues of $A$ are a conjugate pair $\lambda,\bar\lambda$ with $|\lambda|>1$, and $\bar\lambda/\lambda$ is not a root of unity.  By [HP] (or by  Proposition 3.1), the degree sequence $d_1^{(n)} = {\rm deg}(f^n)$ does not satisfy any linear recurrence relation with constant coefficients. On the other hand the inverse map is a polynomial map
$$g:=f_A^{-1} (x_1,x_2,x_3) = (x_3, x_1 x_3,x_2 x_3).$$ 
By Theorem 2.7, the degree sequence $d_k^{-}:= {\rm deg}(f^{-k}_A)$ satisfies $d^-_k = d^-_{k-3}+d^{-}_{k-2}+d^{-}_{k-1}$. 

Finally, we show that in fact $g$ can be made 1-regular in the sense of (1.2).  We start with the induced map on ${\bf P}^3$ 
$$g: [x_0:x_1:x_2:x_3] \mapsto [x_0^2:x_0 x_3:x_1 x_3:x_2 x_3].$$ 
The indeterminacy locus is $ {\it I}(g)=\{x_0=x_3=0\}\cup \{e_3:=[0:0:0:1]\}$.  The orbits of the exceptional hypersurfaces are: $$\eqalign{g : & \{ x_0=0\} \mapsto \{x_0=x_1=0\} \mapsto e_3 \in {\it I}(g)\cr & \{x_3=0\} \mapsto [1:0:0:0] \mapsto [1:0:0:0].}$$  
We see that $g$ does not satisfy (1.2)  by the criterion of [FS]:  an exceptional hypersurface is mapped, after two iterates, completely inside the indeterminacy locus.  We consider the complex manifold $\pi:X\to {\bf P}^3$ obtained by blowing up the point $e_3$ and then the line $\Sigma_{01}= \{x_0=x_1=0\}$.  By $E_3$ and $S_{01}$ we denote the exceptional fibers over $e_3$ and $\Sigma_{01}$ respectively. Under the induced map $g_X$, we have $$g_X:\Sigma_0:= \{x_0 = 0\} \mapsto S_{01} \cap \Sigma_0 \mapsto E_3 \cap \Sigma_{01} \mapsto E_3 \cap \Sigma_0 \mapsto S_{01} \cap \Sigma_0.$$ It follows from [BK, Theorem 1.4] that the induced map $g_X$ is $1$-regular.  And by duality, $g^{-1}_X=\tilde f_A$ is 2-regular.

\bigskip\noindent\centerline{\bf References}
\medskip

\item{[BK]} E. Bedford and KH. Kim, Degree growth of matrix inversion: birational maps of symmetric, cyclic matrices,  arXiv:math/0512507v1 

\item{[BV]} M.P. Bellon and C.-M. Viallet, Algebraic entropy, Comm. Math. Phys.. 204 (1999), 425-437.

\item{[DF]}  J. Diller and C. Favre, Dynamics of bimeromorphic maps of
surfaces, Amer. J. of Math., 123 (2001), 1135--1169.

\item{[F]} C. Favre, Les applications monomiales en deux dimensions, Mich. Math. J.  51 (2003), 467-475.

\item{[FJ]} C. Favre and M. Jonsson, Dynamical compactifications of ${\bf C}^2$.

\item{[FS]} J.-E. Forn\ae ss and N. Sibony, Complex dynamics in higher dimension, II.  Modern Methods in  Complex Analysis (Princeton, NJ, 1992)  Ann. of Math. Stud. Vol. 137, Princeton Univ. Press, Princeton, NJ, 1995, 135--182.

\item{[HP]}  B. Hasselblatt and J. Propp, Degree-growth of monomial maps, arxiv.org/math.DS/0604521

\item{[RS]} A. Russakovskii and B. Shiffman, Value distribution for sequences of rational mappings and complex dynamics, Indiana U. Math.\ J., 46 (1997), 897--932.

\item{[S]}  R.P. Stanley, {\sl Enumerative Combinatorics}, Vol 1. Wadsworth and Brooks/Cole, Monterey, CA, 1986

\bigskip
\rightline{Indiana University}

\rightline{Bloomington, IN 47405}

\rightline{\tt bedford@indiana.edu}

\bigskip
\rightline{Florida State University}

\rightline{Tallahassee, FL 32306}

\rightline{\tt kim@math.fsu.edu}

\bye